\RequirePackage{ifpdf}
\ifpdf % We are running pdfTeX in pdf mode
\documentclass[pdftex]{sigma}
\else
\documentclass{sigma}
\fi

\usepackage{color}

\begin{document}

\allowdisplaybreaks

\renewcommand{\PaperNumber}{035}

\FirstPageHeading

\newtheorem{Theorem}{Theorem}[section]
\newtheorem{Proposition}[Theorem]{Proposition}
\newtheorem{Lemma}[Theorem]{Lemma}
\newtheorem{Open}[Theorem]{Open Question}
\newtheorem{Corollary}[Theorem]{Corollary}
\newtheorem{Conjecture}[Theorem]{Conjecture}
\newtheorem{Specialthm}{Theorem}
\newtheorem{Recipe}[Theorem]{Recipe}

\theoremstyle{definition}
\newtheorem{Question}{Question}
\newtheorem{Motivation}{Motivation}
\newtheorem{Remark}[Theorem]{Remark}
\newtheorem{Definition}[Theorem]{Definition}

\renewcommand\theSpecialthm{\Alph{Specialthm}}

\newcommand{\wt}{\operatorname{wt}}
\newcommand{\bz}{\Bbb{Z}}
\newcommand{\asl}{\widehat{\mathfrak{{sl}}}}
\newcommand{\hook}{\text{hook}}
\newcommand{\arm}{\text{arm}}
\newcommand{\high}{\text{high}}

\renewcommand{\theenumi}{\roman{enumi}}
\renewcommand{\labelenumi}{$(\theenumi)$}

\ShortArticleName{Monomial Crystals and Partition Crystals}

\ArticleName{Monomial Crystals and Partition Crystals}

\Author{Peter TINGLEY}

\AuthorNameForHeading{P.~Tingley}

\Address{Department of Mathematics, Massachusetts Institute of Technology, \\
77 Massachusetts Avenue, Cambridge, MA 02139, USA}
\Email{\href{mailto:ptingley@math.mit.edu}{ptingley@math.mit.edu}}
\URLaddress{\url{http://www-math.mit.edu/~ptingley/}}

\ArticleDates{Received February 10, 2010, in f\/inal form April 12, 2010;  Published online April 21, 2010}

\Abstract{Recently Fayers introduced a large family of combinatorial realizations of the fundamental crystal $B(\Lambda_0)$ for $\asl_n$, where the vertices are indexed by certain partitions. He showed that special cases of this construction agree with the Misra--Miwa realization and with Berg's ladder crystal. Here we show that another special case is naturally isomorphic to a realization using Nakajima's monomial crystal.}

\Keywords{crystal basis; partition; af\/f\/ine Kac--Moody algebra}

\Classification{17B37; 05E10}

\section{Introduction}

Fix $n \geq 3$ and let $B(\Lambda_0)$ be the crystal corresponding to the fundamental representation of $\asl_n$.
Recently Fayers \cite{Fayers:2009} constructed an uncountable family of combinatorial realizations of $B(\Lambda_0)$, all of whose underlying sets are indexed by certain partitions. Most of these are new, although two special cases have previously been studied. One is the well known Misra--Miwa realization~\cite{MM:1990}. The other is the ladder crystal developed by Berg \cite{Berg:2009}.

The monomial crystal was introduced by Nakajima in \cite[Section~3]{Nakajima:2003} (see also~\cite{HN:2006, Kashiwara:2001}). Na\-kaji\-ma considers a symmetrizable Kac--Moody algebra whose Dynkin diagram has no odd cycles, and constructs combinatorial realizations for the crystals of all integrable highest weight modules. In the case of the fundamental crystal $B(\Lambda_0)$ for $\asl_n$, we shall see that the construction works exactly as stated in all cases, including $n$ odd when there is an odd cycle.

Here we construct an isomorphism between a realization of $B(\Lambda_0)$ using Nakajima's monomial crystal and one case of Fayers' partition crystal. Of course any two realizations of $B(\Lambda_0)$ are isomorphic, so the purpose is not to show that the two realizations are isomorphic, but rather to give a simple and natural description of that isomorphism.

This article is organized as follows. Sections~\ref{crystals},~\ref{mono-section} and~\ref{Fayers-section} review necessary background material. Section~\ref{the-iso} contains the statement and proof of our main result. In Section~\ref{questions} we brief\/ly discuss some questions arising from this work.

\section{Crystals} \label{crystals}

In Sections~\ref{mono-section} and~\ref{Fayers-section} we review the construction of Nakajima's monomial crystals and Fayers' partition crystals. We will not assume the reader has any prior knowledge of these constructions. We will however assume that the reader is familiar with the notion of a crystal, so will only provide enough of an introduction to that subject to f\/ix conventions, and refer the reader to  \cite{Kashiwara:1995} or \cite{Hong&Kang:2000} for more details.

We only consider crystals for the af\/f\/ine Kac--Moody algebra $\asl_n$. For us, an $\asl_n$ crystal is the crystal associated to an integrable highest weight $\asl_n$ module. It consists of a set $B$ along with operators $\tilde e_{\bar \imath}, \tilde f_{\bar \imath}: B \rightarrow B \cup \{ 0 \}$ for each $\bar \imath$ modulo $n$, which satisfy various axioms. Often the def\/inition of a crystal includes three functions $\wt, \varphi, \varepsilon: B \rightarrow P$, where $P$ is the weight lattice. In the case of crystals of integrable modules, these functions can be recovered (up shifting in a~null direction) from the knowledge of the $\tilde e_{\bar \imath}$ and $\tilde f_{\bar \imath}$, so we will not count them as part of the data.

\section{The monomial crystal} \label{mono-section}

This construction was f\/irst introduced in \cite[Section 3]{Nakajima:2003}, where it is presented for symmetrizable Kac--Moody algebras where the Dynkin diagram has no odd cycles. In particular, it only works for $\asl_n$ when $n$ is even. However, in Section~\ref{the-iso} we show that for the fundamental crystal $B(\Lambda_0)$ the most naive generalization to the case of odd $n$ gives rise to the desired crystal, so the results in this note hold for all  $n \geq 3$.

We now f\/ix some notation, largely following \cite[Section~3]{Nakajima:2003}.

%$\bullet$ Fix an integer $n \geq 3$.

$\bullet$ Let $\widetilde I$ be the set of pairs $(\bar \imath, k)$ where $\bar \imath$ is a residue mod $n$ and $k \in \bz$.

$\bullet$ Def\/ine commuting formal variables $Y_{\bar \imath,k}$ for all pairs $(\bar \imath, k) \in \widetilde I$.

$\bullet$ Let $\mathcal{M}$ be the set of monomials in the variables $Y_{\bar \imath, k}^{\pm 1}$. To be precise, a monomial $m \in \mathcal{M}$ is a product $\prod\limits_{(\bar \imath, k) \in \widetilde I} Y_{\bar \imath, k}^{u_{\bar \imath, k}}$ with all $u_{\bar \imath, k} \in \bz$ and $u_{\bar \imath, k}=0$ for all but f\/initely many $(\bar \imath, k) \in \widetilde I$.

$\bullet$ For each pair $(\bar \imath, k) \in \widetilde I$, let $A_{\bar \imath,k}= Y_{\bar \imath,k-1} Y_{\bar \imath, k+1} Y_{\bar \imath+\bar 1, k}^{-1} Y_{\bar \imath-\bar 1, k}^{-1}$.

$\bullet$ Fix a monomial $m = \prod\limits_{(\bar \imath, n) \in \widetilde I} Y_{\bar \imath, k}^{u_{\bar \imath, k}} \in \mathcal{M}$. For $L \in \bz$ and $\bar \imath \in \bz / n \bz$, def\/ine:
\begin{gather*}
\wt(m):=\sum_{(\bar \imath, k) \in \widetilde I} u_{\bar \imath, k} \Lambda_{\bar \imath},
\qquad
\varepsilon_{\bar \imath,L}(m) := -{\sum_{l\geq L}}u_{\bar \imath,l}(m),
\qquad
\varphi_{\bar \imath,L}(m) := \sum_{l\leq L} u_{\bar \imath,l}(m),
\\
\varepsilon_{\bar \imath}(m) := \max\{ \varepsilon_{\bar \imath,L}(m) \mid L\in\bz\},
\qquad
p_{\bar \imath}(m) := \max\{L\in\bz\mid \varepsilon_{\bar \imath,L}(m)=\varepsilon_{\bar \imath}(m)\},
\\
\varphi_{\bar \imath}(m) := \max\{ \varphi_{\bar \imath,L}(m) \mid L\in\bz\},
\qquad
q_{\bar \imath}(m) := \min\{L\in\bz\mid \varphi_{\bar \imath,L}(m)=\varphi_{\bar \imath}(m)\}.
\end{gather*}
Note that one always has $\varphi_{\bar \imath}(m), \varepsilon_{\bar \imath}(m) \geq 0$. Furthermore, if $\varepsilon_{\bar \imath}(m)>0$ then $p_{\bar \imath}$ is f\/inite, and if $\varphi_{\bar \imath}(m)>0$ then $q_{\bar \imath}(m)$ is f\/inite.

$\bullet$ Def\/ine $\tilde{e}_{\bar \imath}^M,\tilde{f}_{\bar \imath}^M\colon \mathcal{M} \cup \{ 0 \} \to \mathcal{M} \cup\{0\}$ for each residue $\bar \imath$ modulo $n$
by $\tilde{e}_{\bar \imath}^M(0)= \tilde{f}_{\bar \imath}^M(0)=0$,
\begin{equation} \label{ef-def}
%\begin{split}
   \tilde{e}^M_{\bar \imath}(m) :=
  \begin{cases}
    0 & \text{if $\varepsilon_{\bar \imath}(m) = 0$},
    \\
   A_{\bar \imath,p_{\bar \imath}(m)-1}m & \text{if $\varepsilon_{\bar \imath}(m) > 0$},
  \end{cases}
\qquad \tilde{f}^M_{\bar \imath}(m) :=
  \begin{cases}
    0 & \text{if $\varphi_{\bar \imath}(m) = 0$},
    \\
   A_{\bar \imath,q_{\bar \imath}(m)+1}^{-1}m & \text{if $\varphi_{\bar \imath}(m) > 0$}.
  \end{cases}
%\end{split}
\end{equation}

\begin{Definition} A monomial $m$ is called \emph{dominant} if $u_{\bar \imath, k} \geq 0$ for all $(\bar \imath, k) \in \widetilde I.$
\end{Definition}

\begin{Definition} \label{compatibe-def} Assume $n$ is even. Then
$m$ is called \emph{compatible} if $u_{\bar \imath, k} \neq 0$ implies $k \cong \bar \imath$ modulo~$2$.
\end{Definition}

\begin{Definition}
Let $\mathcal{M}(m)$ be the set of monomials in $\mathcal{M}$ which can be reached from $m$ by applying various $\tilde e_{\bar \imath}^M$ and $\tilde f_{\bar \imath}^M$.
\end{Definition}

\begin{Theorem}[\protect{\cite[Theorem 3.1]{Nakajima:2003}}]  \label{mono-th}
Assume $n >2$ is even, and let $m$ be a dominant, compa\-tib\-le monomial. Then $\mathcal{M}(m)$ along with the operators $\tilde e_{\bar \imath}^M$ and $\tilde f_{\bar \imath}^M$ is isomorphic to the crystal $B({\wt(m)})$ of the integrable highest weight $\asl_n$ module $V({\wt(m)})$.
\end{Theorem}

\begin{Remark}
Notice that although Theorem~\ref{mono-th} only holds when $n$ is even, the operators $\tilde e_{\bar \imath}^M$ and $\tilde f_{\bar \imath}^M$ are well def\/ined for any $n \geq 3$ and any monomial $m$. When $n$ is odd or $m$ does not satisfy the conditions of Theorem~\ref{mono-th}, $\mathcal{M}(m)$ need not be a crystal. However, as we prove in Section~\ref{the-iso}, even when $n$ is odd $\mathcal{M}(Y_{\bar 0, 0})$ is a copy of the crystal $B(\Lambda_0)$ of the fundamental representation of $\asl_n$.
\end{Remark}

\begin{figure}[t]
\centerline{\includegraphics{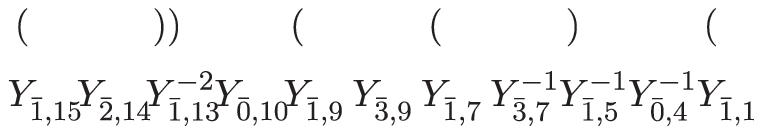}}

\caption{The operators $\tilde e_{\bar \imath}^M$ and $\tilde f_{\bar \imath}^M$ on a monomial $m \in \mathcal{M}$ for $n=4$. We calculate $\widetilde e_{\bar 1}^M$ and $\widetilde f_{\bar 1}^M$. The factors $Y_{\bar \imath, k}$ of $m$ are arranged from left to right by decreasing $k$. The string of brackets $S_{\bar 1}(m)$ is as shown above the monomial. The f\/irst uncanceled ``$)$'' from the right corresponds to a factor of $Y_{\bar 1,13}^{-1}$. Thus $\widetilde e^M_{\bar 1}(m) = A_{\bar 1, 12} m = Y_{\bar 1, 11}^{} Y_{\bar 1, 13}^{}  Y_{\bar 0, 12}^{-1}  Y_{\bar 2, 12}^{-1} m$. The f\/irst uncanceled ``$($'' from the left corresponds to a~factor of $Y_{\bar 1, 9},$ so $\widetilde f^M_{\bar 1}(m) = A_{\bar 1, 10}^{-1} m = Y_{\bar 1, 9}^{-1} Y_{\bar 1, 11}^{-1}  Y_{\bar 0, 10}  Y_{\bar 2, 10} m$.\label{mono-rule}}
\end{figure}

We f\/ind it convenient to use the following slightly dif\/ferent but equivalent def\/inition of $\tilde e_{\bar \imath}^M$ and $\tilde f_{\bar \imath}^M$.
For each $\bar \imath$ modulo $n$, let $S_{\bar \imath}(m)$ be the string of brackets which contains a ``$($" for every factor of $Y_{\bar \imath, k}$ in $m$ and a ``$)$'' for every factor of $Y_{\bar \imath, k}^{-1} \in m$, for all $k \in \bz$. These are ordered from left to right in decreasing order of $k$, as shown in Fig.~\ref{mono-rule}. Cancel brackets according to usual conventions, and set
\begin{gather} \label{ef-def2}
\begin{split}
  & \tilde{e}^M_{\bar \imath}(m)  \hspace{-0.07cm} = \hspace{-0.07cm}
  \begin{cases}
    0 & \mbox{} \hspace{-0.15cm}  \text{if there is no uncanceled ``)'' in $m$},
    \\
   A_{\bar \imath,k-1}m & \mbox{} \hspace{-0.15cm}  \text{if the f\/irst  uncanceled ``)'' from the right comes from a factor $Y_{\bar \imath, k}^{-1}$},
  \end{cases}
\\
  & \tilde{f}^M_{\bar \imath}(m)  \hspace{-0.07cm} = \hspace{-0.07cm}
  \begin{cases}
    0 &  \mbox{} \hspace{-0.15cm}  \text{if there is no uncanceled ``('' in $m$},
    \\
   A_{\bar \imath, k+1}^{-1}m & \mbox{} \hspace{-0.15cm}  \text{if the f\/irst uncanceled ``('' from the left comes from a factor $Y_{\bar \imath, k}$}.
  \end{cases}
\end{split}\!\!\!\!\!\!\!
\end{gather}
It is a straightforward exercise to see that the operators def\/ined in \eqref{ef-def2} agree with those in~\eqref{ef-def}.

\section{Fayers' crystal structures} \label{Fayers-section}

\subsection{Partitions} \label{Partitions}

\begin{figure}[t]
\centerline{\includegraphics{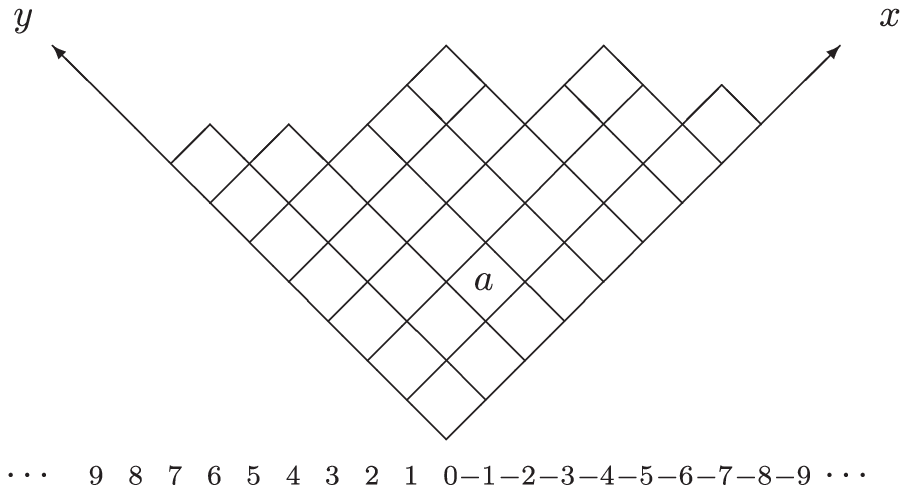}}

\caption{The partition $\lambda=(7,6,5,5,5,3,3,1)$, drawn in ``Russian" notation. The parts $\lambda_i$ are the lengths of the ``rows" of boxes sloping southeast to northwest. The center of each box $b$ has coordinate $(x_b,y_b)$ for some $x_b,y_b \in \bz+1/2$. For the box labeled $a$, $x_a=2.5$ and $y_a=1.5$. The content $c(a)=y_a-x_a$ records the horizontal position of $a$, reading right to left. In this case, $c(a)=-1$. Other relevant statistics are $\hook(a)= 8$, $\arm(a)=3$ and $h(a)= 4$. This partition is not in $B^H$ for $\asl_4$, since the box $a$ is $A^H$-illegal.\label{partition}}
\end{figure}

A partition $\lambda$ is a f\/inite length non-increasing sequence of positive integers. Associated to a~partition is its Ferrers diagram. We draw these diagrams as in Fig.~\ref{partition} so that,
if $\lambda = (\lambda_1,\ldots, \lambda_N)$, then $\lambda_i$ is the number of boxes in row $i$
(rows run southeast to northwest).  Let $\mathcal{P}$ denote the set of all partitions. For $\lambda, \mu \in \mathcal{P}$, we say $\lambda$ is contained in $\mu$ if the diagram
for $\lambda$ f\/its inside the diagram for $\mu$.

Fix $\lambda \in \mathcal{P}$ and a box $b$ in (the diagram of) $\lambda$. We now def\/ine several statistics of $b$. See Fig.~\ref{partition} for an example illustrating these. The \emph{coordinates} of $b$ are the coordinates $(x_b,y_b)$ of the center of $b$, using the axes shown in Fig.~\ref{partition}. The \emph{content} $c(b)$ is $y_b-x_b$. The \emph{arm length} of $b$ is
$\arm(b):= \lambda_{x_b+1/2}-y_b-1/2,$
where $\lambda_{x_b+1/2}$ is the length of the row through $b$. The \emph{hook length} of $b$ is
$\hook(b):= \lambda_{x_b+1/2}-{y_b}+\lambda'_{y_b+1/2}-{x_b},$
where $\lambda_{x_b+1/2}$ is the length of the row containing~$b$ and~$\lambda'_{y_b+1/2}$ is the length of the column containing~$b$.

\subsection{The general construction} \label{pc-section}

We now recall Fayers' construction~\cite{Fayers:2009} of the crystal $B_{\Lambda_0}$ for $\asl_n$ in its most general version. We begin with some notation. An \emph{arm sequence} is a sequence $A=A_1,A_2, \ldots$ of integers such that
\begin{enumerate}\itemsep=0pt
\item $t-1 \leq A_t \leq (n-1) t$ for all $t \geq 1$, and

\item $A_{t+u} \in [A_t+A_u, A_t+A_u+1]$ for all $t,u \geq 1$.
\end{enumerate}
Fix an arm sequence $A$. A box $b$ in a partition $\lambda$ is called \emph{$A$-illegal} if, for some $t \in \bz_{>0}$, $\hook(b)=nt$ and $\arm(b)=A_t$. A partition $\lambda$ is called \emph{$A$-regular} if it has no $A$-illegal boxes. Let $B^A$ denote the set of $A$-regular partitions.

For $\lambda \in \mathcal{P}$ and a box $b \in \lambda$. The \emph{color} of $b$ is the residue $\overline{c}(b)$ of $c(b)$ modulo $n$, where as in Section~\ref{Partitions}, $c(b)$ is the content of $b$. See Fig.~\ref{flat-crystal}.

Fix $\lambda \in \mathcal{P}$ and a residue  $\bar \imath$ modulo $n$. Def\/ine
\begin{itemize}\itemsep=0pt
\item$A(\lambda)$ is the set of boxes $b$ which can be added to $\lambda$ so that the result is still a partition.
\item $R(\lambda)$ is the set of boxes $b$ which can be removed from $\lambda$ so that the result is still a partition.
\item $A_{\bar \imath}(\lambda) = \{ b \in A(\lambda) \text{ such that } \overline c(b)=\bar \imath \}$.
\item $R_{\bar \imath}(\lambda)= \{ b \in R(\lambda) \text{ such that } \overline c(b) = \bar \imath \}$.
\end{itemize}

For each partition $\lambda$ and each arm sequence $A$, def\/ine a total order $\succ_A$ on $ A_{\bar \imath}(\lambda) \cup R_{\bar \imath}(\lambda)$ as follows. Fix $b=(x,y),b'=(x',y') \in A_{\bar \imath}(\lambda) \cup R_{\bar \imath}(\lambda)$, and assume $b \neq b'$. Then there is some $t \in \bz \backslash \{ 0 \}$ such $c(b')-c(b)=nt$. Interchanging $b$ and $b'$ if necessary, assume $t > 0$. Def\/ine $b' \succ_A b$ if $y'-y > A_t$, and $b \succ_A b'$ otherwise.
It follows from the def\/inition of an arm sequence that $\succ_A$ is transitive.

Fix a partition $\lambda$. Construct a string of brackets $S^A_{\bar \imath}(\lambda)$ by placing a ``$($" for every $b \in A_{\bar \imath}(\lambda)$ and a ``$)$" for every $b \in R_{\bar \imath}(\lambda)$, in decreasing order from left to right according to $\succ_A$. Cancel brackets according to the usual rule. Def\/ine maps $\tilde e^A_{\bar \imath}, \tilde f^A_{\bar \imath}: \mathcal{P} \cup \{ 0 \} \rightarrow \mathcal{P} \cup \{ 0 \}$ by $e^A_{\bar \imath}(0)=f^A_{\bar \imath}(0)=0$,
\begin{gather} \label{Hcryst}
\begin{aligned}
& \tilde e_{\bar \imath}^A(\lambda) =
\begin{cases}
\lambda \backslash b~~ & \text{if the f\/irst uncanceled ``$)$'' from the right in $S^A_{\bar \imath}$ corresponds to $b$}, \\
0 & \text{if there is no uncanceled ``$)$'' in $S^A_{\bar \imath}$},
\end{cases} \\
& \tilde f_{\bar \imath}^A(\lambda)=
\begin{cases}
\lambda \sqcup b & \text{if the f\/irst uncanceled ``$($'' from the left in $S^A_{\bar \imath}$ corresponds to $b$}, \\
0 & \text{if there is no uncanceled ``$($'' in $S^A_{\bar \imath}$}.
\end{cases}
\end{aligned}
\end{gather}

\begin{Theorem}[\protect{\cite[Theorem 2.2]{Fayers:2009}}] \label{Fayers-th} Fix $n \geq 3$ and an arm sequence $A$. Then $B^A \cup \{ 0 \}$ is preserved by the maps $\tilde e_{\bar \imath}^A$ and $\tilde f_{\bar \imath}^A$, and forms a copy of the crystal $B(\Lambda_0)$ for $\asl_n$, where as above~$B^A$ is the set of partitions with no $A$-illegal hooks.
\end{Theorem}

\begin{Remark} The operators $\tilde e_{\bar \imath}$ and $\tilde f_{\bar \imath}$ are def\/ined on all partitions. However, as noted in~\cite{Fayers:2009}, the component generated by a non $A$-regular partition need not be an $\asl_n$ crystal.
\end{Remark}

\subsection{Special case: the horizontal crystal} \label{flat-section}
Consider the case of the construction given in Section \ref{pc-section} where, for all $t$, $A_t=\lceil nt/2 \rceil-1$ (it is straightforward to see that this satisf\/ies the def\/inition of an arm sequence). This arm sequence will be denoted $A^H$. For convenience, we denote $B^{A^H}$ simply by $B^H$ and the operators $\tilde e_{\bar \imath}^{A^H}$ and~$\tilde f_{\bar \imath}^{A^H}$ from Section~\ref{pc-section} by $\tilde e_{\bar \imath}^{H}$ and $\tilde f_{\bar \imath}^{H}$.

\begin{Definition} Let $b=(x,y)$ be a box.
The \emph{height} of $b$ is $h(b):= x+y$.
\end{Definition}

\begin{Lemma} \label{flat-order} Fix $\lambda \in B^{H}$, and let $b, b' \in A_{\bar \imath}(\lambda) \cup R_{\bar \imath}(\lambda)$ with $b \neq b'$. Then $b' \succ_{A^H} b$ if and only~if
\begin{enumerate}\itemsep=0pt
\item $h(b') > h(b)$, or
\item  $h(b')= h(b)$ and $c(b') > c(b).$
\end{enumerate}
\end{Lemma}

\begin{proof} Since $\bar c(b) = \bar c(b')$, we have $c(b')-c(b)=nt$ for some $t \in \bz \backslash \{ 0 \}$. First consider the case when $t>0$. Then by def\/inition, $b' \succ_{A^H} b$ if and only if
\begin{equation} \label{newe}
y'-y> A^H_t,
\end{equation}
where
\begin{gather*}
A^H_t  =  \lceil nt/2 \rceil-1 =  \lceil (c(b')-c(b))/2 \rceil-1 = \lceil (y'-x'-y+x)/2 \rceil-1 \\
\phantom{A^H_t}{} = \begin{cases}
(y'-x'-y+x)/2 -1&  \text{if} \ y'-x'-y+x \ \text{is even}, \\
(y'-x'-y+x)/2 -1/2 &  \text{if} \ y'-x'-y+x \ \text{is odd}.
\end{cases}
\end{gather*}
Since $y'-y \in \bz$, \eqref{newe} holds if and only if
$y'-y > (y'-x'-y+x)/2 -1/2$, which rearranges to $h(b') \geq h(b)$.

The case $t<0$ follows immediately from the case $t>0$ since $\succ_{A^H}$ is a total order.
\end{proof}

\begin{figure}[t]
\centerline{\includegraphics{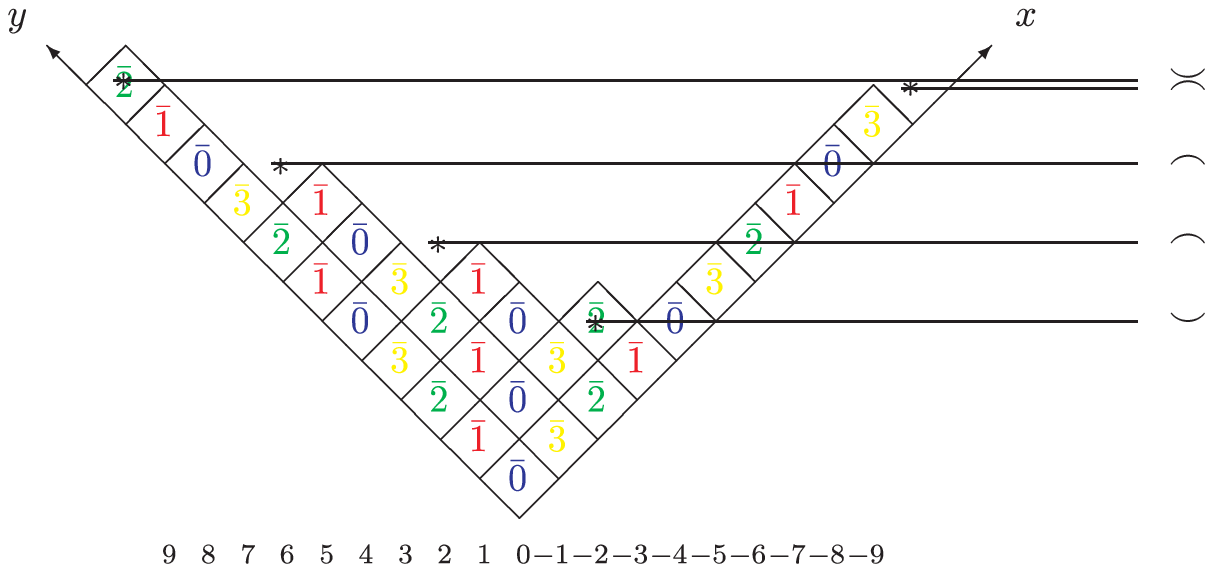}}

\caption{The partition $\lambda=(11,7,4,2,1,1,1,1,1,1)$ is an element of $B^H$ for $\asl_4$, since no box $b$ in~$\lambda$ is $A^H$-illegal. The color $\overline c(b)$ of each box $b$ is written inside $b$. We demonstrate the calculation of~$\tilde f_{\bar 2}(\lambda)$. The string of brackets $S_{\bar 2}^H(\lambda)$ has a ``('' for each {\color{green} ${\bar 2}$}-addable box and a ``)'' for each ${\color{green} {\bar 2}}$-removable box, ordered from left to right lexicographically, f\/irst by decreasing height $h(b)$, then by decreasing content~$c(b)$. The result is shown on the right of the diagram (rotate the page 90 degrees counter clockwise). Thus $\tilde f_{\bar 2}(\lambda)$ is the partition obtained by adding the box with coordinates $(x,y)= (10.5,0.5)$. The map $\Psi$ from Section~\ref{the-iso} takes $\lambda$ to
$Y_{\bar 3, 11}^{} Y_{\bar 2, 8}^{}  Y_{\bar 2,6}^{}  Y_{\bar 3,5}^{}  Y_{\bar 1, 5}^{}  Y_{\bar 2, 10}^{}  Y_{\bar 2, 12}^{-1} Y_{\bar 1, 9}^{-1} Y_{\bar 1, 7}^{-1} Y_{\bar 2, 6}^{-1} Y_{\bar 3, 11}^{-1}$. After reordering and simplifying, this becomes $Y_{\bar 2, 12}^{-1}Y_{\bar 2, 10}^{} Y_{\bar 1, 9}^{-1} Y_{\bar 2, 8}^{}  Y_{\bar 1, 7}^{-1} Y_{\bar 1, 5}^{} Y_{\bar 3, 5}^{}.$ The string of brackets $S^M_{\bar 2}(\Psi(\lambda))$ is the same as the string of brackets $S^H_{\bar 2} (\lambda)$, except that a canceling pair $()$ has been removed. The condition that $\lambda$ has no $A^H$-illegal boxes implies that $S_{\bar \imath}^H(\lambda)$ and $S_{\bar \imath}^M(\Psi(\lambda))$ are always  the same up to removing pairs of canceling brackets, which is essentially the proof that $\Psi$ is an isomorphism.\label{flat-crystal} }
\end{figure}

Lemma \ref{flat-order} implies that $\tilde e_{\bar \imath}^{H}$ and $\tilde f_{\bar \imath}^{H}$ are as described as in Fig.~\ref{flat-crystal}.

\section{A crystal isomorphism} \label{the-iso}

Def\/ine a map $\Psi: \mathcal{P} \cup \{ 0 \} \rightarrow \mathcal M \cup \{ 0 \}$ by $\Psi(0)=0$, and, for all $\lambda \in \mathcal{P}$,
\begin{equation*}
\Psi(\lambda) := \prod_{b \in A(\lambda)} Y_{\bar c(b), h(b)-1}  \prod_{b \in R(\lambda)} Y_{\bar c(b), h(b)+1}^{-1}.
\end{equation*}
Here $A(\lambda)$ and $R(\lambda)$ are as in Section~\ref{pc-section}.

\begin{Theorem} \label{iso-th}
For any $n \geq 3$, any $\bar \imath$ modulo $n$, and any $\lambda \in B^H$, we have $\Psi(\tilde e_{\bar \imath}^H \lambda) = \tilde e_{\bar \imath}^M \Psi(\lambda)$ and $\Psi(\tilde f_{\bar \imath}^H \lambda) = \tilde f_{\bar \imath}^M \Psi(\lambda)$, where $\tilde e_{\bar \imath}^M$, $\tilde f_{\bar \imath}^M$ are as in Section~{\rm \ref{mono-section}} and $\tilde e_{\bar \imath}^H$, $\tilde f_{\bar \imath}^H$ are as in Section~{\rm \ref{flat-section}}.
\end{Theorem}

Before proving Theorem~\ref{iso-th} we will need a few technical lemmas.

\begin{Lemma} \label{bv}
Let $\lambda$ and $\mu$ be partitions such that $\mu = \lambda \sqcup b$ for some box $b$. Then $\Psi(\mu) = A_{\overline c(b), h(b)}^{-1} \Psi(\lambda).$
\end{Lemma}

\begin{proof} Let $i= \bar c(b)$.
It is clear that the pair $( A(\lambda), R(\lambda))$ dif\/fers from  $( A(\mu), R(\mu))$ in exactly the following four ways:
\begin{itemize}\itemsep=0pt
\item $b \in A_{\bar \imath}(\lambda) \backslash A_{\bar \imath}(\mu)$.

\item   $b \in R_{\bar \imath}(\mu) \backslash R_{\bar \imath}(\lambda)$.

\item Either $(i)$: $A_{\bar \imath + \bar 1}(\mu) \backslash A_{\bar \imath + \bar 1}(\lambda)= b'$ and $ R_{\bar \imath + \bar 1}(\lambda) =R_{\bar \imath + \bar 1}(\mu)$ for some box $b'$ with $h(b') = h(b)+1$, or $(ii)$: $R_{\bar \imath + \bar 1}(\lambda) \backslash R_{\bar \imath + \bar 1}(\mu) = b'$ and $ A_{\bar \imath + \bar 1}(\lambda) =A_{\bar \imath + \bar 1}(\mu)$ for some box $b'$ with $h(b')= h(b)-1$.

\item Either $(i)$: $A_{\bar \imath - \bar 1}(\mu) \backslash A_{\bar \imath - \bar 1}(\lambda)= b''$ and $ R_{\bar \imath - \bar 1}(\lambda) =R_{\bar \imath - \bar 1}(\mu)$ for some box $b''$ with $h(b'') = h(b)+1$, or $(ii)$: $R_{\bar \imath - \bar 1}(\lambda) \backslash R_{\bar \imath - \bar 1}(\mu) = b''$ and $ A_{\bar \imath - \bar 1}(\lambda) =A_{\bar \imath - \bar 1}(\mu)$ for some box $b''$ with $h(b'')= h(b)-1$.
\end{itemize}
By the def\/inition of $\Psi$, this implies
\begin{gather*}
\Psi(\mu)= Y_{\bar \imath, h(b)-1}^{-1}  Y_{\bar \imath, h(b)+1}^{-1}  Y_{\bar \imath+\bar 1, h(b)}^{}  Y_{\bar \imath-\bar 1, h(b)}^{}  \Psi(\lambda) = A_{\bar \imath, h(b)}^{-1} \Psi(\lambda).\tag*{\qed}
\end{gather*}\renewcommand{\qed}{}
\end{proof}

\begin{Lemma} \label{when-illegal}
Let $\lambda \in B^H$, and choose $b \in A_{\bar \imath}(\lambda)$, $b' \in R_{\bar \imath}(\lambda)$. Then
\begin{enumerate}\itemsep=0pt
\item \label{oneoff} $h(b) \neq h(b')+1$.

\item \label{twooff} If $h(b)= h(b')$ then $c(b') >c(b)$, so $b' \succ_{A^H} b$.

%\item \label{threeoff} If $h(b)= h(b')+2$ then $c(b) >c(b')$.

\item \label{neworder} $b \succ_{A^H} b'$ if and only if $h(b)-1 \geq h(b')+1$.

\item \label{no-others} Assume $h(b)-1=h(b')+1$. If $c \in A_{\bar \imath}(\lambda)$ satisfies $b \succ_{A^H} c \succ_{A^H} b'$, then $h(c)=h(b)$.

\item \label{no-others2} Assume $h(b)-1=h(b')+1$. If $c \in R_{\bar \imath}(\lambda)$ satisfies $b \succ_{A^H} c \succ_{A^H} b'$, then $h(c)=h(b')$.
\end{enumerate}
\end{Lemma}

\begin{figure}[t]
\centerline{\includegraphics{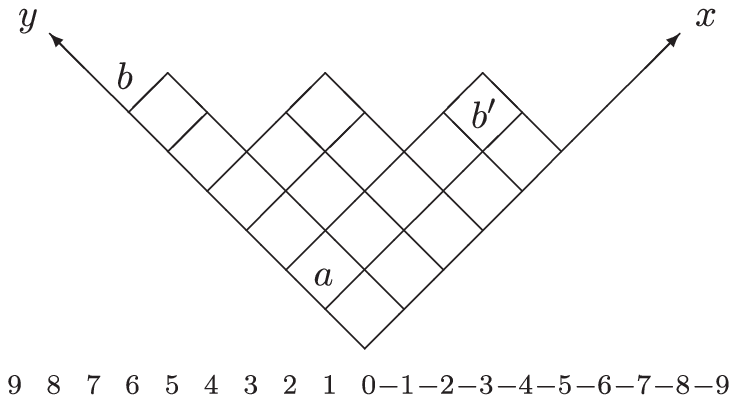}}

\caption{The hook def\/ined by $b \in A_{\bar \imath}(\lambda)$ and $b' \in R_{\bar \imath}(\lambda)$. There will always be a unique box $a$ in $\lambda$ which is in either the same row or the same column as $b$ and also in either the same row or the same column as~$b'$. Taking $n=3$, we have $b \in A_{\bar 0}(\lambda)$ and $b' \in R_{\bar 0}(\lambda)$. Then $\hook(a) =9 =3 \times 3$ and $\arm(a)=4 = A^H_3$, so $a$ is $A^H$-illegal. It is straightforward to see that in general, if either $(i)$:~$h(b)=h(b')+1$ or $(ii)$:~$h(b)=h(b')$ and $c(b)>c(b')$, then the resulting hook is $A^H$-illegal, and hence $\lambda \notin B^H$. \label{hook-fig}}
\end{figure}

\begin{proof}
By the def\/initions of $A_{\bar \imath}(\lambda)$ and $R_{\bar \imath}(\lambda)$, $b$ and $b'$ cannot lie in either the same row or the same column, which implies that there is a unique box $a$ in $\lambda$ which shares a row or column with each of $b$, $b'$. It is straightforward to see that  if~$(\ref{oneoff})$ or $(\ref{twooff})$ is violated then this $a$ is $A^H$-illegal (see Fig.~\ref{hook-fig}).

To see part $(\ref{neworder})$, recall that by Lemma~\ref{flat-order}, $b \succ_{A^H} b'$ if and only if $h(b) > h(b')$ or both $h(b)=h(b')$ and $c(b)> c(b')$. This order agrees with the formula in part $(\ref{neworder})$, since parts $(\ref{oneoff})$ and $(\ref{twooff})$ eliminate all cases where they would dif\/fer.

Part $(\ref{no-others})$ and $(\ref{no-others2})$ follow because any other $c \in A_{\bar \imath}(\lambda) \cup R_{\bar \imath}(\lambda)$ with $b \succ_{A^H} c \succ_{A^H} b'$ would violate either $(\ref{oneoff})$ or $(\ref{twooff})$.
\end{proof}

\begin{proof}[Proof of Theorem~\ref{iso-th}]
Fix $\lambda \in B^H$ and $\bar \imath \in \bz/n\bz$. Let $S^M_{\bar \imath}(m)$ denote the string of brackets used in Section~\ref{mono-section} to calculate $\tilde e_{\bar \imath}^M(m)$ and $\tilde f_{\bar \imath}^M(m)$.
Let $S_{\bar \imath}^H(\lambda)$ denote the string of brackets used in Section~\ref{Fayers-section} to calculate  $\tilde e_{\bar \imath}^H(\lambda)$ and $\tilde f_{\bar \imath}^H(\lambda)$, and def\/ine the height of a bracket in $S^H_{\bar \imath}(\lambda)$ to be~$h(b)$ for the corresponding box $b \in A_{\bar \imath}(\lambda) \cup R_{\bar \imath}(\lambda)$.

By Lemma \ref{when-illegal} parts $(\ref{no-others})$ and $(\ref{no-others2})$, for each $k \geq 1$, all ``$($''  in $S^H_{\bar \imath}(\lambda)$ of height $k+1$ are immediately to the left of all ``$)$'' of height $k-1$.  Let $T$ be the string of brackets obtained from $S^H_{\bar \imath}(\lambda)$ by, for each $k$, canceling as many ``$($'' of height $k+1$ with ``$)$" of height $k-1$ as possible. Notice that one can use $T$ instead of $S^H_{\bar \imath}(\lambda)$ to calculate $\tilde e_{\bar \imath}^H(\lambda)$ and $\tilde f_{\bar \imath}^H(\lambda)$ without changing the result.

By the def\/inition of $\Psi$, it is clear that
\begin{enumerate}\itemsep=0pt
\item The ``$($'' in $T$ of height $k+1$ correspond exactly to the factors of $Y_{\bar \imath, k}$ in $\Psi(\lambda)$.

\item The ``$)$'' in $T$ of height $k-1$ correspond exactly to the factors of $Y_{\bar \imath, k}^{-1}$ in $\Psi(\lambda)$.
\end{enumerate}

Thus the brackets in $T$ correspond exactly to the brackets in $S^M_{\bar \imath}(\Psi(\lambda))$. Furthermore, Lemma~\ref{when-illegal} part~$(\ref{neworder})$ implies that these brackets occur in the same order. The theorem then follows from Lemma~\ref{bv} and the def\/initions of the operators (see equations~\eqref{ef-def2} and~\eqref{Hcryst}).
\end{proof}

\begin{Corollary}
For any $n$, $\mathcal{M}(Y_{\bar 0,0})$ is a copy of the fundamental crystal $B(\Lambda(0))$ for $\asl_n$.
\end{Corollary}

\begin{proof}
This follows immediately from Theorem \ref{iso-th}, since, by Theorem \ref{Fayers-th}, $B^H$ is a copy of the crystal $B(\Lambda_0)$.
\end{proof}

\section{Questions} \label{questions}

\begin{Question} Nakajima originally developed the monomial crystal using the theory of $q$-cha\-rac\-ters from~\cite{FR:1999}. Can this theory be modif\/ied to give rise to any of Fayers' other crystal structures? One may hope that this would help explain algebraically why these crystal structures exist.
\end{Question}

\begin{Question}
In \cite{Kim:2005}, Kim considers a modif\/ication to the monomial crystal developed by Kashiwara \cite{Kashiwara:2001}. She works with more general integral highest weight crystals, but restricting her results to $B(\Lambda_0)$ one f\/inds a natural isomorphisms between this modif\/ied monomial crystal and the Misra--Miwa realization. The Misra--Miwa realization corresponds to one case of Fayer's partition crystal, but not the one studied in Section~\ref{flat-section}. In \cite{Kashiwara:2001}, there is some choice as to how the monomial crystal is modif\/ied. Do other modif\/ications also correspond to instances of Fayers crystal? Which instances of Fayers' partiton crystal correspond to modif\/ied monomial crystals (or appropriate generalizations)?
\end{Question}

\begin{Question}
The monomial crystal construction works for higher level~$\asl_n$ crystals. There are also natural realizations of higher level $\asl_n$ crystals using tuples of partitions (see~\cite{FLOTW:1997, JMMO:1991, KangLee, abacus}). Is there an analogue of Fayers' construction in higher levels generalizing both of these types of realization?
\end{Question}

\subsection*{Acknowledgments}
We thank Chris Berg, Matthew Fayers, David Hernandez and Monica Vazirani for interesting discussions. This work was supported by NSF grant DMS-0902649.

\pdfbookmark[1]{References}{ref}
\LastPageEnding


\begin{thebibliography}{99}

\footnotesize\itemsep=0pt

\bibitem{Berg:2009}
Berg C.,
$(\ell,0)$-JM partitions and a ladder based model for the basic  crystal of $\widehat{\mathfrak{sl}}_\ell$,
\href{http://arxiv.org/abs/0901.3565}{arXiv:0901.3565}.

\bibitem{Fayers:2009}
Fayers M.,
Partition models for the crystal of the basic ${U}_q(\widehat{\frak{sl}}_n)$-module,
\href{http://dx.doi.org/10.1007/s10801-010-0217-9}{{\it J. Algebraic Combin.}}, to appear,
\href{http://arxiv.org/abs/0906.4129}{arXiv:0906.4129}.

\bibitem{FLOTW:1997}
Foda O.,, Leclerc B., Okado M., Thibon J.-Y., Welsh~T.A.,
Branching functions of {$A\sp {(1)}\sb {n-1}$} and Jantzen--Seitz problem for Ariki--Koike algebras,
\href{http://dx.doi.org/10.1006/aima.1998.1783}{{\em Adv. Math.}} {\bf 141} (1999), 322--365,
\href{http://arxiv.org/abs/q-alg/9710007}{q-alg/9710007}.

\bibitem{FR:1999}
Frenkel E., Reshetikhin N.,
The $q$-characters of representations of quantum af\/f\/ine algebras and deformations of $\mathcal W$-algebras,
in Recent Developments in Quantum Af\/f\/ine Algebras and Related Topics (Raleigh, NC, 1998),  {\em Contemp. Math.}, Vol.~248, Amer. Math. Soc., Providence, RI, 1999, 163--205,
\href{http://arxiv.org/abs/math.QA/9810055}{math.QA/9810055}.

\bibitem{HN:2006}
Hernandez D., Nakajima H.,
Level 0 monomial crystals,
\href{http://projecteuclid.org/euclid.nmj/1167159343}{{\em Nagoya Math. J.}} {\bf 184} (2006), 85--153,
\href{http://arxiv.org/abs/math.QA/0606174}{math.QA/0606174}.

\bibitem{Hong&Kang:2000}
Hong J., Kang S.-J.,
Introduction to quantum groups and crystal bases,
{\em Graduate Studies in Mathematics}, Vol.~42, Amer. Math. Soc., Providence, RI, 2002.

\bibitem{JMMO:1991}
Jimbo M., Misra K.C., Miwa T., Okado M.,
Combinatorics of representations of $\text{U}_q(\widehat{\frak{sl}}(n))$ at $q=0$,
\href{http://dx.doi.org/10.1007/BF02099073}{{\em Comm. Math. Phys.}} {\bf 136} (1991), 543--566.

\bibitem{KangLee}
Kang S.-J., Lee H.,
Higher level af\/f\/ine crystals and Young walls,
\href{http://dx.doi.org/10.1007/s10468-006-9013-6}{{\em Algebr. Represent. Theory}} {\bf 9} (2006), 593--632,
 \href{http://arxiv.org/abs/math.QA/0310430}{math.QA/0310430}.

\bibitem{Kashiwara:1995}
Kashiwara M.,
On crystal bases,
 in Representations of Groups (Banf\/f, AB, 1994), {\it CMS Conf. Proc.}, Vol.~16, Amer. Math. Soc., Providence, RI, 1995, 155--197.

\bibitem{Kashiwara:2001}
Kashiwara M.,
Realizations of crystals,
in  Combinatorial and Geometric Representation Theory (Seoul, 2001),   {\em Contemp. Math.}, Vol.~325,  Amer. Math. Soc., Providence, RI, 2003, 133--139,
\href{http://arxiv.org/abs/math.QA/0202268}{math.QA/0202268}.

\bibitem{Kim:2005}
Kim J.-A.,
Monomial realization of crystal graphs for ${U}_q({A}_n(1))$,
\href{http://dx.doi.org/10.1007/s00208-004-0613-3}{{\em Math. Ann.}} {\bf 332} (2005), 17--35.

\bibitem{MM:1990}
Misra K., Miwa T.,
Crystal base for the basic representation of $U_q(\widehat{\mathfrak{sl}}_n)$,
\href{http://dx.doi.org/10.1007/BF02102090}{{\em Comm. Math. Phys.}} {\bf 134} (1990), 79--88.

\bibitem{Nakajima:2003}
Nakajima H.,
$t$-analogs of $q$-characters of quantum af\/f\/ine algebras of type $A\sb n$, $D\sb n$,
in  Combinatorial and Geometric Representation Theory (Seoul, 2001),   {\em Contemp. Math.}, Vol.~325,  Amer. Math. Soc., Providence, RI, 2003, 141--160,
\href{http://arxiv.org/abs/math.QA/0204184}{math.QA/0204184}.

\bibitem{abacus}
Tingley P.,
Three combinatorial models for $\widehat{\rm sl}_n$ crystals, with applications to cylindric plane partitions,
\href{http://dx.doi.org/10.1093/imrn/rnm143}{{\em Int. Math. Res. Not. IMRN}} {\bf 2008} (2008), no.~2, Art.~ID rnm143, 40~pages,
\href{http://arxiv.org/abs/math.QA/0702062}{math.QA/0702062}.

\end{thebibliography}
\end{document}